\documentclass[a4paper,11pt]{article}

\usepackage{amsfonts} 
\usepackage{amssymb} 
\usepackage{amsmath,amsthm} 

\newtheorem{theoreme}{Theorem}[section] 
\newtheorem{lemma}[theoreme]{Lemma} 
\newtheorem{prop}[theoreme]{Proposition} 
\newtheorem{cor}[theoreme]{Corollary} 
\newtheorem{definition}[theoreme]{Definition}

\def \Rk {\ {\bf Remark.} } 
\def \Rks {\ {\bf Remarks.}} 
\def \sm {\setminus } 
\newcommand{\be}{\begin{enumerate}}  \newcommand{\ee}{\end{enumerate}} 
\newcommand{\bi}{\begin{itemize}}  \newcommand{\ei}{\end{itemize}} 
\newcommand{\bd}{\begin{description}}  \newcommand{\ed}{\end{description}}

\newcommand{\comment}[1]{}

\def \R {\mathbb{R}}  \def \N {\mathbb{N}}

\DeclareMathOperator{\vol}{vol}
\DeclareMathOperator{\Cone}{Cone}

\numberwithin{equation}{section}  
\renewcommand{\phi}{\varphi} 
\renewcommand{\epsilon}{\varepsilon} 
 
\title{On manifolds with quadratic curvature decay}
\author{Nader Yeganefar}
\date{\today} 

\begin{document} 
\maketitle
\begin{abstract}
We give conditions which imply that a complete noncompact manifold with quadratic curvature decay has finite topological type. In particular, we find links between the topology of a manifold with quadractic curvature decay and some properties of the asymptotic cones of such a manifold.
\end{abstract}

\section{Introduction}
Let $M$ be a complete noncompact Riemmannian manifold. In this paper, we are interested in the following basic problem: find geometric conditions which imply that $M$ has finite topological type, i.e. is homeomorphic to the interior of a compact manifold with boundary. For example, it is known that if $M$ is flat or has nonnegative sectional curvature, then it has finite topological type. (In fact much stronger results on the structure of $M$ hold in this case by the work of Cheeger and Gromoll \cite{CG}.) 

In studying our basic problem, an important general principle to keep in mind is that what is really relevant should be the geometry at infinity. For instance, a manifold which is flat or nonnegatively curved only outside a compact subset has finite topological by \cite{GW}. We might then expect that this conclusion also holds for manifolds which are in some sense "asymptotically flat" or "asymptotically nonnegatively curved". Indeed, Abresch \cite{A} generalized the results of Greene and Wu in the following form. Assume that for some $m_0$ in $M$ and some constants $C,\epsilon >0$, the sectional curvatures $K$ of $M$ at all points $m$ satisfy
$$K\geq -C/d(m_0,m)^{2+\epsilon},$$
where $d$ denotes the distance function on $M$. Then $M$ has finite topological type. Moreover, this theorem is optimal in the sense that on any noncompact surface Abresch constructed a complete metric for which we have 
$$K=o(1/d(m_0,m)^2),$$
as $m$ goes to infinity.

We now introduce the class of manifolds which will be the focus of this paper.
\begin{definition}
Let $(M,m_0)$ be a pointed complete noncompact Riemannian manifold. 
\bi
\item[1)] We say that $M$ has lower quadratic curvature decay if for some $C>0$ the sectional curvatures $K$ of $M$ at all points $m$ satisfy
$$K\geq -C/d(m_0,m)^2.$$
\item[2)] We say that $M$ has quadratic curvature decay if for some $C>0$ the sectional curvatures $K$ of $M$ at all points $m$ satisfy
$$\vert K\vert\leq C/d(m_0,m)^2.$$
\ei
\end{definition}
\noindent \Rks
\bi
\item[i)] If $(M,m_0)$ has (lower) $C$-quadratic curvature decay and $m'_0$ is any point in $M$, then it is easy to see that $(M,m'_0)$ has (lower) $C'$-quadratic curvature decay for some constant $C'>0$.
\item[ii)] Having (lower) $C$-quadratic curvature decay is independent of constant rescalings of the metric.
\item[iii)] J. Lott \cite{L} has a slightly different definition of quadratic curvature decay. Namely he considers the condition
$$\limsup _{r\to\infty} {\sup _{m, d(m_0,m)=r}{r^2\vert K\vert}}\leq C.$$
Qualitatively, this is the same definition as ours, but quantitatively it is a bit more general. 
\ei
The examples of Abresch mentioned above show that having quadratic curvature decay does not restrict the toplogy of surfaces. We quote also this striking result of Gromov (see \cite[Lemma 2.1]{LS}): on \emph{any} noncompact manifold there exists a complete metric of quadratic curvature decay. Therefore we need to find additional assumptions that restrict the topology of a manifold carrying a metric of quadratic curvature decay. This is done for example by Lott-Shen \cite{LS}, Sha-Shen \cite{SS}, Lott \cite{L}, do Carmo-Xia \cite{dCX}, Xia \cite{X}, etc.

Before stating our main technical result, let us first say a few words about asymptotic cones; more details will be given in the next section. If $(M,m_0,d_M)$ is a pointed metric space and if $\{ R_i\}$ is a sequence of positive numbers going to infinity, we can consider the sequence of rescaled pointed metric spaces $\{(M,m_0,d_M/R_i)\}$. If this sequence is precompact in the pointed Gromov-Hausdorff topology, then any of its limit points is called an asymptotic cone. Intuitively, an asymptotic cone has to reflect the large scale metric behaviour of $M$. Moreover, an asymptotic cone may not be unique (i.e may depend on the converging subsequence of the original sequence) and may even not be a metric cone, see \cite{CC}, \cite{M2}. Now, even if $\{(M,m_0,d_M/R_i)\}$ is not precompact, there is a contruction using "ultrafilters" and "ultralimits" which allows us to get from this sequence a pointed metric space denoted by $\Cone _{\omega ,\{R_i\}}{(M,m_0)}$, where $\omega$ is a nonprincipal ultrafilter. We call this space also asymptotic cone, and it indeed generalizes the first definition given above. Recall finally that a metric space  $M$ is said to have a pole at some point $m_0$ if for any point $m$ there exists a ray (i.e. a minimizing geodesic on $[0,\infty )$) starting at $m_0$ and passing through $m$. 
We can now state our first result
\begin{theoreme}\label{main}
Let $M$ be complete noncompact manifold with lower quadratic curvature decay. If $M$ has infinite topological type, then there is a sequence of positive numbers $\{R_i\}$ diverging to infinity such that $\Cone _{\omega , \{R_i\}}{(M,m_0)}$ does not have a pole at its basepoint.
\end{theoreme}
This result may be seen as a generalization of the first part of the main theorem of Petrunin and Tuschmann \cite{PT}. Let us give a rough idea of the proof. Recall that there is a notion for a point in $M$ to be critical for the (not smooth) distance function $d(m_0,.)$. This is due to Grove and Shiohama (see the surveys of Cheeger \cite{C} and Grove \cite{Gr}). For a manifold with lower quadratic curvature decay, we will prove a distance estimate at each critical point (Lemma \ref{fundamental}). This estimate is the analog of the one obtained by Sormani in another context, see the "uniform cut lemma" \cite[Lemma 7]{So2}. If we assume that our manifold has infinite topological type, then it follows by critical point theory that there is an infinite sequence of critical points $q_i$ such that $R_i=d(m_0,q_i)$ is diverging. Then we prove our theorem by slightly modifying the proof of \cite[Theorem 11]{So2}.
 
In order to have some applications of our theorem, it would be interesting to know large classes of manifolds for which all asymptotic cones have a pole. We consider here the case of noncompact complete manifolds of nonnegative Ricci curvature. It is well-known that these manifolds have at least linear volume growth and at most Euclidean volume growth of geodesic balls. 
%More precisely, for $m_0$ in $M$ and $R>0$, denote by $(B(m_0,R)$ the geodesic ball of radius $R$ and center $R$, and by $\vol (B(m_0,R))$ its volume. Then for some constant $c>0$ we have
%$$cR\leq \vol (B(m_0,R))\leq \omega_nR^n.$
%Here $n$ is the dimension of $M$ and $\omega _n$ is the volume of the unit sphere in $\R^n$. The first inequality is due to Calabi-Yau \cite{Y} and the second one is a consequence of
Moreover, it was shown by Cheeger-Colding \cite{CG} and Sormani \cite{So1} that if $M$ has respectively Euclidean volume growth or linear volume growth, then every asymptotic cone is a metric cone, and hence has a pole at its basepoint. Therefore, we get from Theorem \ref{main} the following corollary (whose first case was obtained by Sha-Shen \cite{SS} more directly):
\begin{cor}
Let $M$ be a complete Riemannian manifold with basepoint $m_0$. Assume that  $M$ has lower quadratic curvature decay and nonnegative Ricci curvature. If
 $M$ has either Euclidean volume growth or linear volume growth, then $M$ has finite topological type.
\end{cor}
\noindent\Rks
\bi
\item[i)] There are examples of manifolds of nonnegative Ricci curvature, Euclidean volume growth, and infinite topological type, see \cite{M1}
\item[ii)] As already mentioned, the case of Euclidean volume growth was obtained also by Sha and Shen \cite[Theorem 1.1]{SS}. Sha and Shen also treat a case which is close to minimal volume growth \cite[Theorem 1.2]{SS}. Namely, for $r>0$ set
$$v(r)=\inf _{m\in B(m_0,r)}\vol (B(m,1)),$$
where for $B(m,r)$ denotes the ball of radius $r$ centered at $m$ and $\vol (B(m,r))$ denotes its volume. Sha and Shen showed that if $M$ has nonnegative Ricci curvature and satisfies
\begin{equation}\label{ShaShen}
\limsup _{r\to \infty} \frac{\vol (B(m_0,r))}{v(r)r^2}=0,
\end{equation} 
then it has finite topologocal type (the condition on Ricci curvature is actually not necessary here, as was shown later by Lott and Shen \cite{LS}). It seems hard to have a control on $v(r)$ even on a manifold with nonnegative Ricci curvature so that it is not clear that (\ref{ShaShen}) implies linear volume growth.
\item[iii)] Sha and Shen asked if it is true that every manifold with lower quadratic curvature decay and nonnegative Ricci curvature has finite topological type. Actually this doesn't hold and counterexamples were constructed by Menguy. 
\ei
Next, we apply the technics used to prove Theorem \ref{main} to get (implicit) estimates of the criticality radius on some manifolds. Recall that the criticality radius at some point $m_0$ is the largest $R\in(0,\infty]$ such that there is no critical point of $d(m_0,.)$ (other than $m_0$) in the geodesic ball $B(m_0,R)$. It is always bigger than the injectivity radius at that point. 
\begin{theoreme}\label{cr}
Given constants $n\in \N$, $C, \Lambda, v>0$, there exists $R=R(n,C,\Lambda,v)>0$ with the following property: if $(M,m_0)$ is a pointed $n$-dimensional Riemannian manifold  such that
\be
\item $M$ has $C$-quadratic curvature decay,
\item $M$ has bounded sectional curvature $\vert K \sb M \vert \leq \Lambda \sp 2$,
\item $M$ has Euclidean volume growth:  $\forall m\in M, \forall t\in \R ,\,\vol {(B(m,t))} \geq vt^n$,
\ee
then the criticality radius at $m_0$ is greater than or equal to R.
Moreover, for fixed $n$, $\Lambda$ and $v$, the function $R$ goes to infinity when $C$ goes to zero.
\end{theoreme}
Note that our assumptions 1 and 3 above are scale invariant, so that we cannot hope to get similar results without 2 or at least without some extra assumption which is not scale invariant. We would like also to emphasize here that the existence of $R$ is actually already known. Namely, on the one hand Klingenberg \cite{Kli} gave a general lower bound on the injectivity radius $i$ at $m_0$ in the following form
$$i\geq \min{(l,\pi /\Lambda )},$$
where $l$ denotes the length of the shortest geodesic loop at $m_0$. For this estimate, no lower bound on the volume growth is needed, nor any quadratic curvature decay assumption. On the other hand, Cheeger, Gromov and Taylor \cite{CGT} gave a lower bound on $l$ in terms of $\Lambda$, a lower bound on the volume growth, and the dimension of the manifold. These together give a lower bound on $i$ under the assumptions $2$ and $3$ of our theorem, and hence a lower bound on the criticality radius at $m_0$. However, this lower bound is at most $\pi /\Lambda$ and under the assumptions of $C-$quadratic curvature decay and $C$ small, we get a qualitative improvement of it. 

The organization of the paper is as follows. In the next next, we recall the necessary background material on ultralimits and asymptotic cones. In Section 3, we prove our distance estimate for critical points and deduce from it a slightly more general form of Theorem \ref{main}. In the last section, we prove first the existence of the function $R$ in Theorem \ref{cr} (actually under less restrictive assumptions) and then finish the proof of Theorem \ref{cr}.

\noindent {\bf Acknowledgements.} I would like to thank Gilles Carron, Guofang Wei and John Lott for useful comments.

%%%%%%%%%%%%%%%%%%%%%%%%%%%%%%%%%%%%%%%%%%%%%%%%%%%%%%%%%%%%%%%%%
%%%%%%%%%%%%%%%%%%%%%%%%%%%%%%%%%%%%%%%%%%%%%%%%%%%%%%%%%%%%%%%%%
%%%%%%%%%%%%%%%%%%%%%%%%%%%%%%%%%%%%%%%%%%%%%%%%%%%%%%%%%%%%%%%%%%

\section{Asymptotic cones}

\subsection{Ultralimits}
In this section, we recall standard facts about ultralimits. The material is taken from \cite[Chapter 9]{K} and \cite[Section 2.4]{KlL} and the reader should consult these sources for further references and developments.

A nonprincipal ultrafilter is a finitely additive probability measure $\omega$ on the subsets of $\N$ such that for every $I\subset \N$, we have $\omega (I)=0$ or $1$, and $\omega (I)=0$ if $I$ is finite. It follows easily from the definition that if two subets $I$ and $J$ have full measure, then $I\cap J$ has also full measure.

If $Y$ is a compact metric space and if $f : \N \to Y$ is a map, then we can take the "limit" of $f$ with respect to $\omega$; this is the unique element $y\in Y$, denoted by $\omega -\lim {f}$, such that for every neighborhood $U$ of $y$, the preimage $f^{-1}(U)$ has full $\omega$-measure.

Consider now a sequence $(M_i,m_i,d_{M_i})$ of pointed metric spaces. Any ultrafilter $\omega$ allows us to put these spaces together and get a pointed metric space $(M_\omega ,m_\omega ,d_{M_\omega })$, which is called the \emph{ultralimit} of the sequence. To define this object, we proceed as follows. Let $M_\infty$ be the space of sequences $\{p_i\}_{i\in \N}$, with $p_i \in M_i$, such that the sequence $\{d_{M_i}(m_i, p_i)\}$ is bounded. For two elements $\{p_i\}$ and $\{q_i\}$ of $M_\infty$, we can in particular consider the number $d_\infty(\{p_i\},\{q_i\})$ defined by
$$d_\infty(\{p_i\},\{q_i\}) = \omega -\lim {\{ i\mapsto d_{M_i}(p_i, q_i)\} }.$$
$d_\infty$ is a pseudodistance on $M_\infty$ and we define
$$(M_\omega ,d_{M_\omega})=(M_\infty ,d_\infty)/\sim \, ,$$
where we identify two elements of $M_\infty$ whose $d_\infty$-distance is zero. Finally, the sequence $\{m_i\}$ defines an element of $M_\omega$ which we denote by $m_\omega$.

The following proposition shows that ultralimits are generalizations of Gromov-Hausdorff limits (see \cite[Propositions 9.2 and 9.4]{K} and \cite[Lemma 2.4.3]{KlL} for the proof, and \cite{G} for background material on Gromov-Hausdorff convergence)
\begin{prop}[\cite{K}, \cite{KlL}]\label{ultra}
Let $(M_i,m_i,d_{M_i})$ be a sequence of pointed metric spaces and let $\omega$ be a nonprincipal ultrafilter on $\N$.
%\be
%\item If each $M_i$ is a geodesic metric space, then $M_\omega$ is also a geodesic space.
%\item 
If $(M_i,m_i,d_{M_i})$ is a precompact family in the pointed Gromov-Hausdorff topology, then the metric space $(M_\omega, m_\omega, d_{M_\omega})$ is a limit point of this family in the pointed Gromov-Hausdorff topology.
%\ee
\end{prop}

%%%%%%%%%%%%%%%%%%% Asymptotic cones  %%%%%%%%%%%%%%%

\subsection{Asymptotic cones}

Ultralimits are particularly useful to define the notion of \emph{asymptotic cone} for a metric space.
\begin{definition}\label{cone}
Let $(M,m_0,d_M)$ be a pointed metric space and let $\{R_i\}$ be a sequence of positive real numbers diverging to infinity. For each $i$, let $(M_i,m_0,d_M/R_i)$ be the space $(M,m_0)$ with the rescaled metric $d_M/R_i$. The generalized asymptotic cone of $(M,m_0)$, with respect to $\{R_i\}$ and a given nonprincipal ultrafilter $\omega$, is defined by
$$\Cone _{\omega ,\{R_i\}}{(M,m_0)}=\omega -\lim {(M,m_0,d_M/R_i)}.$$
\end{definition}
\noindent \Rks
\bi
\item[i)] We use the word "generalized" in the definition to avoid confusion with the more traditional notion of asymptotic cone, as explained in the introduction. Namely, assume that the sequence $(M,m_0,d_M/R_i)$ is precompact in the pointed Gromov-Hausdorff topology. (By the Gromov compactness theorem, this is for example the case if $M$ has nonnegative Ricci curvature \cite[Theorem 5.3]{G}.) Then there is a subsequence of this sequence which converges to a metric space $(X,x,d_X)$, called also the asymptotic cone of $(M,m_0)$. This space $X$ may not be a metric cone (see \cite{M2}) and may even not be unique (i.e. may depend on the convergent subsequence, see \cite{CC}). However, by Proposition \ref{ultra} above, for a suitable convergent subsequence, $X$ will be the same as the asymptotic cone introduced in Definition \ref{cone}.
\item[ii)] The usual choice of the sequence $R_i$ is $R_i=i$, but this is not necessary. The definition of $\Cone _{\omega ,\{R_i\}}{(M,m_0)}$ actually makes sense for any positive sequence of real numbers $\{R_i\}$, and not only for divergent sequences; nevertheless, for a general sequence, the terminology "asymptotic cone" is probably not a good one. Furthermore, if we have a sequence $(M_i,m_i,d_{M_i})$ of metric spaces, we can define its (generalized) asymptotic cone by
$$\Cone _{\omega ,\{R_i\}}{(M_i,m_i)_{i\in \N}}=\omega -\lim {(M_i,m_i,d_{M_i}/R_i)}.$$
\ei
%It is often an interesting question to know whether $\Cone _{\omega ,\{R_i\}}{(M_i,m_i)_{i\in \N}}$ is a metric cone or not. From the work of Cheeger-Cloding \cite{CC} and our preceding remarks, we deduce the following
%\begin{theoreme}[\cite{CC}]\label{metriccone}
%Let $(M_i^n,m_i)$ be a sequence of $n-$dimensional complete (pointed) Riemannian manifolds with nonnegative Ricci curvature. Let $\{R_i\}$ be a sequence of positive numbers and $\omega$ a nonprincipal ultrafilter on $\N$. Assume that:\\
%either $\{R_i\}$ is bounded and the sequence $M_i$ is volume non collapsing, i.e.
%$$\inf _{i} {\vol {(B_{M_i}(m_i,1))}}>0,$$
%or the sequence $M_i$ has uniform Euclidean volume growth, i.e.
%$$\inf _{i}{\lim _{R\to \infty}R^{-n}{\vol {(B_{M_i}(m_i,R))}}}>0.$$
%Then $\Cone _{\omega ,\{R_i\}}{(M_i,m_i)_{i\in \N}}$ is a metric cone. In particular, it has a pole at its basepoint.
%\end{theorme}
%\begin{proof}
%By the first remark after Definition \ref{cone}, we may assume that \\
%$\Cone _{\omega ,\{R_i\}}{(M_i,m_i)_{i\in \N}}$ is the pointed Gromov-Hausdorff  limit of some subsequence of the rescaled manifolds $(M_j,m_j,d_{M_j}/R_j)$. These manifolds have nonnegative Ricci curvature, and both sets of assumptions of the proposition imply easily that they are volume non collapsing. The conclusion then follows from \cite[Theorem 5.2]{CC}
%\end{proof}

%%%%%%%%%%%%%%%%%%%%%%%%%%%%%%%%%%%%%%%%%%%%%%%%%%%%%%%%%%

\section{A distance estimate for critical points}

To prove Theorem \ref{main}, our first task will be to get a distance estimate for critical points, in the spirit of \cite[Lemma 2.1]{SS} (see also \cite{LS}). For manifolds with nonnegative Ricci curvature, Sormani obtained such an estimate for some special critical points which she called "halfway points" (see \cite{So2} for the explanation of this terminology). Under the assumption of lower quadratic curvature decay, we derive this estimate for all critical points, and not only halfway points.
\begin{lemma}\label{estimate}
Let $M$ be a complete Riemannian manifold with basepoint $m_0$. Assume that for  some $C>0$, $M$ has lower $C$-quadratic curvature decay. Then there exists $\epsilon =\epsilon (C)>0$ with the following property. If $q$ is a critical point for $d(m_0,.)$  and if $m$ is a point with $d(m_0,m)\geq (1+\epsilon)d(m_0,q)$, then
$$d(m,q)> d(m_0,m)-d(m_0,q)+\frac{\epsilon}{2} d(m_0,q).$$
\end{lemma}
\begin{proof}
We argue by contradiction and assume that the lemma is not true. Then for all $\epsilon>0$ we find a critical point $q$ at distance $d(m_0,q)=R$ and a point $m$ at distance $d(m_0,m)=\alpha R$ for some $\alpha \geq 1+\epsilon$, such that 
$$d(m,q)\leq \alpha R-R+\frac{\epsilon}{2}R.$$
We can rewrite this as
\begin{equation}\label{contradiction}
\frac{3}{2}\epsilon R \leq R(\alpha -1+2\epsilon )-d(m,q).
\end{equation}
Now let $\gamma$ be a minimal geodesic from $q$ to $m$. As $q$ is critical for $d(m_0,.)$, there exists a minimal geodesic $\sigma$ from $q$ to $m_0$ such that the angle between $\gamma '(0)$ and $\sigma '(0)$ is $\leq \pi/2$.
Consider the points $x=\sigma (\epsilon R )$ and $y=\gamma (\epsilon R)$. By the triangle inequality, we have
$$d(m_0,m)-d(m_0,x)-d(m,y)\leq d(x,y).$$
Using the fact that $d(m_0,x)=R-\epsilon R$ and $d(m,y)=d(m,q)-\epsilon R$, it follows
\begin{equation}\label{triangle}
R(\alpha -1+2\epsilon )-d(m,q)\leq d(x,y).
\end{equation}
Now inequality \ref{contradiction} implies that if $\epsilon$ is sufficiently small, any minimal geodesic joining $x$ and $y$ is contained in $M\sm B(m_0,R/2)$. Namely, let $c$ be such a geodesic. For each $t$, we have by the triangle inequality
\begin{equation}\label{chaispas}
\alpha R=d(m_0,m)\leq d(m_0,c(t))+d(c(t),y)+d(y,m).
\end{equation}
Moreover, we have
\begin{eqnarray*}
d(c(t),y) &\leq& d(x,y)\\
          &\leq& d(x,q)+d(q,y)\leq 2\epsilon R.   
\end{eqnarray*}
Combining this with inequality \ref{chaispas}, we obtain
$$\alpha R\leq d(m_0,c(t))+2\epsilon R+ d(m,q)-\epsilon R$$
Using our assumption \ref{contradiction} on $d(m,q)$, it follows finally that
$$d(m_0,c(t))\geq(1-\frac{3}{2}\epsilon )R,$$
hence $d(m_0,c(t))$ is bigger than $R/2$ for small $\epsilon$. Moreover, it is  clear that $\gamma \vert _{[0,\epsilon R]}$ and $\sigma \vert _{[0,\epsilon R]}$  are also contained in $M\sm B(m_0,R/2)$. Now we have $K\geq -4C^2/R^2$ on $M\sm B(m_0,R/2)$, so by the Toponogov comparison theorem applied to the hinge $(\gamma \vert _{[0,\epsilon R]}, \sigma \vert _{[0,\epsilon R]})$, inequalities \ref{triangle} and \ref{contradiction}, we get
$$ \cosh (2C\frac{3}{2}\epsilon )\leq \cosh ^2(2C\epsilon ).$$
This is impossible if $\epsilon$ is small enough.
\end{proof}
Now we come to the proof of Theorem \ref{main}. We will use the arguments of \cite[Theorem 11]{So2} with some minor changes. The main technical difference is that Sormani deals with manifolds with nonnegative Ricci curvature which form a precompact family in the pointed Gromov-Hausdorff topology, whereas we work with families which are not necessarily precompact. We have then to use ultralimits. Theorem \ref{main} is a direct consequence of the following result:
\begin{prop}\label{fundamental}
Let $(M_i,m_i)$ be a sequence of complete pointed Riemmannian manifolds with the same lower quadratic curvature decay. Assume that for every $i$, there is a critical point $q_i$ for $d_{M_i}(m_i,.)$ at distance $R_i=d_{M_i}(m_i,q_i)>0$. 
Then for any nonprincipal ultrafilter $\omega$, the space $\Cone _{\omega ,\{R_i\}}{(M_i,m_i)_{i\in \N}}$ does not have a pole at its basepoint.
\end{prop}
\begin{proof}
For simplicity, we set $$X_i=(M_i,m_i,d_{M_i}/R_i),$$ 
$$(M_\omega,m_\omega)= \Cone _{\omega ,\{R_i\}}{(M_i,m_i)_{i\in \N}},$$
 and $q_\omega=[\{ q_i\} ]$. First we note that for all $i$, we have 
$$d_{X_i}(m_i,q_i)=d_{M_i}(m_i,q_i)/R_i=1.$$
As $m_\omega$ is by definition represented by $\{m_i \}$, it follows that $d_{M_\omega}(m_\omega ,q_\omega )=1.$

We will show that there is no ray emanating from $m_\omega$ which passes through $q_\omega$. For this, we argue by contradiction and assume that there is a ray $\gamma : [0,\infty )\to M_\omega$ such that $\gamma (0)=m_\omega$ and $\gamma (1)=q_\omega$. For $\delta >0$, consider the point $p_\omega =\gamma (1+\delta)$, and choose a sequence $p_i\in M_i$ such that $p_\omega =[\{ p_i \} ]$. Let $\epsilon >0$ be as in lemma \ref{estimate} and choose $\eta \in (0, \epsilon /4)$. We have $d_{M_\omega}(m_\omega ,p_\omega)=1+\delta$, because $\gamma$ is minimizing. Therefore, for a set $J=\{j\}$ of indices of full $\omega$-measure, we have
$$d_{X_j}(m_j,p_j)>1+\delta -\eta ,$$
so
$$d_{M_j}(m_j,p_j)>(1+\delta -\eta)R_j.$$
We may assume that $\delta$ and $\eta$ are chosen such that $\delta -\eta \geq \epsilon$. Then by lemma \ref{estimate} (with $M=M_j$, $m_0=m_j$, $q=q_j$ and $m=p_j$ in the notation of this lemma), we get
$$d_{M_j}(p_j,q_j)>(\delta -\eta +\epsilon /2)R_j,$$
which is equivalent to
\begin{equation}\label{contr}
d_{X_j}(p_j,q_j)>\delta -\eta +\epsilon /2.
\end{equation}
On the other hand, we have $d_{M_\omega}(p_\omega ,q_\omega )=\delta$ because $p_\omega =\gamma (1+\delta)$, $q_\omega =\gamma (1)$ and $\gamma$ is minimizing. From this we deduce the existence of a set $K=\{k\}$ of indices of full $\omega$-measure such that
$$d_{X_k}(p_k, q_k)\leq \delta +\eta .$$
Finally, we note that the set $J\cap K$ has still full $\omega$-measure, so that the last inequality combined with \ref{contr} gives
$$\delta -\eta +\epsilon /2<\delta +\eta.$$
This contradicts our choice of $\eta <\epsilon /4$ and finishes the proof of the lemma.
\end{proof}
\Rk The proof shows actually a slightly more general result. Namely, any minimizing geodesic emanating from $m_\omega$ which passes through $q_\omega$ stops being minimizing at a time $t\leq 1+\epsilon$, where $\epsilon$ is as in lemma \ref{estimate}.

%\begin{cor}\label{critic1}
%Given constants $n\in \N$, $C>0$ and $\alpha >0$, there exists $R=R(n,C,\alpha )>0$ with the following property: if $(M,m_0)$ is a pointed $n-$dimensional Riemannian manifold  with lower $C$-quadratic curvature decay, $Ricci \geq 0$, and Euclidean volume growth $\alpha _M\geq \alpha$, then all the critical points of $d(m_0,.)$ are contained in the ball $B(m_0,R)$. 
%\end{cor}
%\begin{proof}
%For the existence of $R$, we argue by contradiction and assume on the contrary that for some constants $n\in \N$, $C>0$ and $\alpha >0$, we can find a sequence $(M_i,m_i)$ of $n$-dimensional manifols such that: 
%\bi
%\item[i)] for all $i$, $M_i$ has $C$-quadratic curvature decay, nonnegative Ricci curvature and $\alpha _{M_i}\geq \alpha$,
%\item[ii)] and for all $i$, there is a critical point of $d_{M_i}(m_i,.)$ at distance $R_i$, with $R_i\to \infty$.
%\ei
%Then by Lemma \ref{fundamental}, the space $\Cone _{\omega ,\{R_i\}}{(M_i,m_i)_{i\in \N}}$ does not have a pole at its basepoint, whereas by Proposition \ref{metriccone}, it has one. This is a contradiction.
%\end{proof}

%%%%%%%%%%%%%%%%%%%%%%%% criticality radius%%%%%%%%%%%%%%%%%%%%

\section{Estimates of the criticality radius}
In this section, we prove Theorem \ref{cr}. The first part of this theorem (i.e. the existence of $R$) is a direct consequence of the following more general result, whose proof is related to some arguments of \cite{L}:
\begin{prop}\label{critic1}
Given constants $n\in \N$, $C, \Lambda, v>0$, there exists $R=R(n,C,\Lambda,v)>0$ with the following property: if $(M,m_0)$ is a pointed $n$-dimensional Riemannian manifold  such that
\be
\item $M$ has lower $C$-quadratic curvature decay,
\item $M$ has bounded sectional curvature $\vert K \sb M \vert \leq \Lambda \sp 2$,
\item $M$ is volume non collapsing (i.e. $\forall m\in M, \vol {(B(m,1))} \geq v$),
\ee
then the criticality radius at $m_0$ is greater than or equal to R.
\end{prop}
\begin{proof}
Assume on the contrary that the proposition is not true. Then there exist constants $n, C, \Lambda, v$ and there exists a sequence of pointed Riemannian $n$-dimensional manifolds $(M_i,m_i,g_{M_i})$ such that for all $i$
\be
\item $M_i$ has lower $C$-quadratic curvature decay,
\item $M_i$ has bounded sectional curvature $\vert K \sb {M_i} \vert \leq \Lambda \sp 2$,
\item $M_i$ is volume non collapsing $\inf _{m\in M_i}{\vol {(B(m,1))}}\geq v$,
\item there exists a critical point $q_i$ for $d_{M_i}(m_i,.)$ at some distance $R_i=d_{M_i}(m_i,q_i)>0$, with $R_i\to 0$.
\ee
Setting $$X_i =(M_i,m_i,g_{M_i}/R_i\sp 2),$$
it follows from Proposition \ref{fundamental} that for any nonprincipal ultrafilter $\omega$ the space $\Cone _{\omega ,\{R_i\}}{(M_i,m_i)_{i\in \N}}$ does not have a pole at its basepoint. We are going to show that this space is in fact isometric to $\R\sp n$ to get a contradiction. 
First, due to the rescaling to define $X_i$, we have
\begin{equation}\label{flat}
\vert K_{X_i}\vert =R_i\sp 2 \vert K_{M_i}\vert \leq R_i\sp 2\Lambda \sp 2.
\end{equation}
In particular, the sequence $\{ X_i\}$ has uniformly bounded sectional curvature because $R_i$ goes to zero. By the Gromov compactness theorem and Proposition \ref{ultra} a subsequence of $(X_i,m_i)$ converges in the pointed Gromov-Hausdorff topology to $\Cone _{\omega ,\{R_i\}}{(M_i,m_i)_{i\in \N}}$.

Moreover, for fixed $R>0$ and $p\in X_i$
$$\vol {(B_{X_i}(p,R))}=\frac{\vol {(B_{M_i}(p,RR_i))}}{R_i\sp n}.$$
We may assume that for $i$ sufficiently large we have $RR_i\leq 1$. If we denote by $V(r)$ the volume of a ball of radius $r$ in the $n$-dimensional hyperbolic space of constant sectional curvature $-\Lambda \sp 2$, the lower bound on the sectional curvature of $M_i$ and the Bishop-Gromov volume comparison theorem imply that
\begin{eqnarray}\label{volume}
\vol {(B_{X_i}(p,R))} &\geq& \frac{\vol {(B_{M_i}(m,1))}}{V(1)} \frac{V(RR_i)}{(RR_i)\sp  n}R\sp n \nonumber \\
                      &\geq& k R\sp n,
\end{eqnarray}
where $k>0$ is a constant depending only on $n$, $\Lambda$ and $v$. Thus we are in the non collapsing case of Gromov-Hausdorff convergence and a subsequence of $(X_i, m_i)$ converges to $\Cone _{\omega ,\{R_i\}}{(M_i,m_i)_{i\in \N}}$ in the pointed $C\sp{1,\beta}$ topology (for any $\beta \in (0,1)$). By \cite{P}, the bounds \ref{flat} on the sectional curvatures and the fact that $R_i$ goes to zero, $\Cone _{\omega ,\{R_i\}}{(M_i,m_i)_{i\in \N}}$ is a flat Riemannian manifold. The volume estimate \ref{volume} and the $C\sp{1,\beta}$ convergence imply that $\Cone _{\omega ,\{R_i\}}{(M_i,m_i)_{i\in \N}}$ has also Euclidean volume growth and is therefore isometric to $\R\sp n$. This finishes the proof.
\end{proof}

\begin{proof}[Proof of Theorem \ref{cr}]
The existence of $R$ is an immediate consequence of Proposition \ref{critic1}. For the last assertion of the theorem, we argue by contradiction as in the proof of Proposition \ref{critic1}. Then there exist constants $n, \Lambda, v$, $R_0>0$ and  a sequence of pointed Riemannian $n$-dimensional manifolds $(M_i,m_i,g_{M_i})$ such that for all $i$
\be
\item $M_i$ has $1/i$-quadratic curvature decay,
\item $M_i$ has bounded sectional curvature $\vert K \sb {M_i} \vert \leq \Lambda \sp 2$,
\item $M_i$ has Euclidean volume growth: $\forall m\in M_i, \forall t\in \R,  \vol {(B_{M_i}(m,t))} \geq vt^n$,
\item there exists a critical point $q_i$ for $d_{M_i}(m_i,.)$ at some distance $R_i=d_{M_i}(m_i,q_i)>0$, with $R_i\leq R_0$.
\ee
We set
 $$X_i =(M_i,m_i,g_{M_i}/R_i\sp 2),$$
so that by Proposition \ref{fundamental} we know that for any nonprincipal ultrafilter $\omega$ the space $\Cone _{\omega ,\{R_i\}}{(M_i,m_i)_{i\in \N}}$ does not have a pole at its basepoint. To derive a contradiction, we will show as in the proof of Proposition \ref{critic1} that this space is isometric to $\R ^n$. First, using the fact that the $M_i$'s have uniformly bounded sectional curvature and the upper bound $R_i\leq R_0$, we deduce that the $X_i$'s have also uniformly bounded sectional curvature. Moreover, it is clear that the $X_i$'s have also uniformly Euclidean volume growth. It follows that a subsequence of $(X_i,m_i)$ converges to $\Cone _{\omega ,\{R_i\}}{(M_i,m_i)_{i\in \N}}$ in the pointed $C^{1,\beta}$ topology (for any $\beta \in (0,1)$). As each $M_i$ has $1/i$-quadratic curvature decay, so does also each $X_i$. Then it is easy to see that $\Cone _{\omega ,\{R_i\}}{(M_i,m_i)_{i\in \N}}$ is a flat $n$-dimensional manifold (see also \cite[Lemma 2]{L}). Furtheremore, $\Cone _{\omega ,\{R_i\}}{(M_i,m_i)_{i\in \N}}$ has also Euclidean volume growth by $C^{1,\beta}$ convergence. Hence it is isometric to $\R ^n$. 
\end{proof}

%%%%%%%%%%%%%%%%%%%%%%%%%%% biblio %%%%%%%%%%%%%%%%%%%

\end{document}